\documentclass[a4paper,12pt]{article}

\usepackage{amsthm}
\usepackage{amssymb,amsmath}
\usepackage{indentfirst}
\usepackage{wasysym}
\usepackage[dvips]{graphicx}

\theoremstyle{plain}
\newtheorem{thm}{Theorem}

\newcommand{\ket}[1]{\vert#1\rangle}
\newcommand{\bra}[1]{\langle#1\vert}
\newcommand{\Tr}{\text{Tr}}
\newcommand{\CCC}{\mathbb{C}}
\newcommand{\RRR}{\mathbb{R}}
\newcommand{\NNN}{\mathbb{N}}
\newcommand{\EEE}{\mathbb{E}}
\newcommand{\PPP}{\mathbb{P}}

\newcommand{\scp}[2]{\langle #1| #2 \rangle}
\newcommand{\Hilbert}{\mathcal{H}}
\newcommand{\sphere}{\mathcal{S}}

\newcommand{\sys}{{\mathrm{sys}}}
\newcommand{\env}{{\mathrm{env}}}

\newcommand{\z}[1]{{#1}}
\newcommand{\y}[1]{{#1}}

\title{Elementary Proof for Asymptotics of Large Haar-Distributed Unitary Matrices}
\author{
Christian Mastrodonato\footnote{Dipartimento di Fisica
dell'Universit\`a di
    Genova and INFN sezione di Genova, Via Dodecaneso 33, 16146
    Genova, Italy.  E-mail: christian.mastrodonato@ge.infn.it} \, and
Roderich Tumulka\footnote{Department of Mathematics,
    Rutgers University, 110 Frelinghuysen Road, Piscataway, NJ 08854-8019,
    USA.  E-mail: tumulka@math.rutgers.edu}
}
\date{October 1, 2007}

\begin{document}

\maketitle

\begin{abstract}
We provide an elementary proof for a theorem due to
Petz and R\'effy which states that for a
random $n\times n$ unitary matrix with distribution given by the
Haar measure on the unitary group $U(n)$, the upper left (or any
other) $k\times k$ submatrix converges in distribution,
after multiplying by a normalization factor
$\sqrt{n}$ and as $n\to\infty$, to a matrix of independent complex Gaussian random
variables with mean 0 and variance 1.

\medskip


\end{abstract}

\section{Introduction}

The aim of this paper is to give an alternative, elementary proof of a theorem first established by Petz and R\'effy in \cite{PR04}, concerning the joint distribution of the upper left $k\times k$ entries of a random unitary $n \times n$ matrix in the limit $n\to\infty$ and formulated as Theorem~1 below. This theorem is of particular interest in quantum statistical mechanics, where \z{one often studies the behavior of a small system (corresponding to dimension $k$) coupled to a heat bath---a much larger system corresponding to dimension $n$. Specifically, Theorem 1} can be used for studying the distribution of the conditional wave function of a system coupled to a heat bath in the relevant limit (in which the size of the heat bath tends to infinity). As we show in \cite{GLMTZ07}, this distribution typically converges, as a consequence of Theorem~1, to the so-called ``GAP'' measure \cite{GLTZ06}, which can thus be regarded as the thermal equilibrium distribution of the conditional wave function. \z{We explain this application further in Section~\ref{sec:GAP}.}

We fix some notation and terminology. Let $\PPP$ denote probability and $\EEE$ expectation, $U(n)$ the group of unitary $n\times n$ matrices, and Haar$(U(n))$ the (normalized) \emph{Haar measure} on this group, representing the ``uniform'' probability distribution over $U(n)$. We write $(a_{ij})$ for the matrix with entries $a_{ij}$. The relevant notion of convergence of probability distributions is \emph{weak convergence}, also known as ``convergence in distribution'' of the random variables \cite[Sec.~25]{Bill}. By a complex Gaussian random variable $G$ with mean $0$ and variance $\sigma^2$ we mean $G=X+iY$, where $X$ and $Y$ are independent real Gaussian random variables with means $\EEE \,X = 0$ and $\EEE \,Y =0$ and variances $\EEE \, X^2 = \sigma^2/2$ and $\EEE \, Y^2 = \sigma^2/2$.

\begin{thm}
If $(U_{ij})$ is Haar$(U(n))$ distributed, then the upper left (or, in fact, any) $k\times k$ submatrix, multiplied by a normalization factor $\sqrt{n}$,
converges in distribution, as $n\to\infty$, to a
random $k\times k$ matrix $(G_{ij})$ whose entries $G_{ij}$ are
independent complex Gaussian random variables with mean 0 and
variance $\EEE |G_{ij}|^2=1$.
\end{thm}

To understand the factor $\sqrt{n}$, note that a column of a unitary $n \times n$ matrix is a unit vector, and thus a single entry should be of order $1/\sqrt{n}$.  A random $k \times k$ matrix such as $(G_{ij})$, consisting of independent complex Gaussian variables with mean 0 and variance 1, is also called ``$\sqrt{k}$ times a standard non-selfadjoint Gaussian matrix.''

Theorem~1 is a generalization of the familiar fact that the first $k$ entries of a random unit vector in $\RRR^n$ (with uniform probability distribution over the unit sphere), multiplied by a normalization factor $\sqrt{n}$, converge in distribution to a vector whose $k$ entries are independent real Gaussian random variables with mean 0 and variance 1.\footnote{As a physical interpretation of this fact, consider $N$ classical particles without interaction in a box $\Lambda\subseteq \RRR^3$; a given energy corresponds to a surface in phase space $\Lambda^N \times \RRR^{3N}$ given by $\Lambda^N \times \sphere$, where $\sphere$ is the sphere of appropriate \y{radius $\propto \sqrt{N}$ in momentum space $\RRR^{3N}$; assuming a random phase point with micro-canonical distribution (i.e., uniform on $\Lambda^N \times \sphere$)}, the marginal distribution of the momentum of the first particle is, \z{in the limit $N\to\infty$,} Gaussian. This fact is part of the justification of \y{Maxwell's law of the Gaussian} distribution of momenta.} This fact \y{(with $\RRR^n$ replaced by $\CCC^n$)} is contained in Theorem~1 by specializing to just the first columns of the matrices $(U_{ij})$ and $(G_{ij})$.

The proof of Petz and R\'effy is based on the convergence of the joint
distribution of the eigenvalues of a $k\times k$ submatrix of an unitary
matrix to the corresponding distribution for a $k\times k$ Gaussian
matrix. Our proof, in contrast, is based on the geometric
properties of Gaussian random matrices. While it involves some more
cumbersome estimates, it employs only elementary methods.

\section{Application to Typicality of GAP Measures}
\label{sec:GAP}

We briefly describe the application of Theorem~1 in quantum statistical mechanics.

Consider a quantum system entangled to its environment, so that the composite has a wave function $\psi \in\Hilbert_\sys \otimes \Hilbert_\env$, with $\Hilbert_\sys$ and $\Hilbert_\env$ the Hilbert spaces of the system and the environment.
Suppose $\Hilbert_\sys$ has dimension $k$, while $\Hilbert_\env$ has very large dimension $n$. According to the Schmidt decomposition, every $\psi \in \Hilbert_\sys \otimes \Hilbert_\env$ can be written as
\begin{equation}
\psi = \sum_{i=1}^k c_i \, \chi_i \otimes \phi_i
\end{equation}
with coefficients $c_i\in\CCC$, an orthonormal basis $\{\chi_1,\ldots,\chi_k\}$ of $\Hilbert_\sys$ and an orthonormal system $\{\phi_1,\ldots,\phi_k\}$ in $\Hilbert_\env$. Relative to any fixed orthonormal basis $\{b_1,\ldots, b_n\}$ of $\Hilbert_\env$, the coefficients $U_{ij} = \scp{b_j}{\phi_i}$ of the $\phi_i$ form the first $k$ rows of an $n\times n$ unitary matrix, and the uniform distribution over all $\psi$'s with a given reduced density matrix
\begin{equation}
\rho_\sys = \sum_i |c_i|^2 \, \ket{\chi_i}\bra{\chi_i}
\end{equation}
gives rise to (\y{the appropriate} marginal of) the Haar measure on $(U_{ij})$.

For reasons we explain below, it is of interest to consider, for a fixed but typical $\psi$, a random column of $(U_{ij})$, or, equivalently, the random vector (arising from a random choice of $j$)
\begin{equation}\label{condwf}
\psi_\sys  = \sum_i c_i \, U_{ij} \, \chi_i =\scp{b_j}{\psi}_\env \in \Hilbert_\sys\,,
\end{equation}
where the scalar product is a partial scalar product. By Theorem~1, in the limit $n\to\infty$, each column of $(U_{ij})$ has a Gaussian distribution, and any two columns are independent; as a consequence, by the law of large numbers, for typical $\psi$ the empirical distribution of $\psi_\sys$ approximates a Gaussian distribution on $\Hilbert_\sys$ with covariance $\rho_\sys$.\footnote{This fact is similar to Maxwell's law in the classical setting of Footnote 1: For a typical phase point on $\Lambda^N \times \sphere$, the empirical distribution of the momenta (over all $N$ particles) approximates a Gaussian distribution on $\RRR^3$ as $N\to\infty$. This follows using the law of large numbers from the fact, described in Footnote 1, that the momentum of each particle is Gaussian-distributed, and that the momenta of different particles are independent.}

This fact is significant for the proof that \emph{the thermal equilibrium distribution of the conditional wave function is the GAP measure}, a particular probability distribution on the unit sphere of Hilbert space. Let us explain.

The notion of \emph{conditional wave function} \cite{DGZ92} is a precise mathematical version of the concept of collapsed wave function.
Conditional on the state $b_j$ of the environment, the conditional wave function $\psi_\sys$ of the system is given by the expression \eqref{condwf} (times a normalizing factor). 
Now replace $j$ by a random variable $J$ with the quantum theoretical probability distribution
\begin{equation}
\PPP(J=j) = \bigl\| \scp{b_j}{\psi}_\env \bigr\|^2\,.
\end{equation}
The resulting random vector $\psi_\sys$ is called the conditional wave function. For example, a system after a quantum measurement is still entangled with the apparatus, but its collapsed wave function is a conditional wave function.

Now consider a system kept in thermal equilibrium at a temperature $1/\beta$ by a coupling to a large heat bath. Even if we assume that $\psi \in \Hilbert_\sys \otimes \Hilbert_\env$ (with the environment being the heat bath) is non-random, the conditional wave function $\psi_\sys$ is random, and for typical $\psi$ within the microcanonical ensemble (i.e., for most $\psi$ relative to the uniform distribution over the subspace corresponding to a narrow energy interval), the distribution of $\psi_\sys$ is a universal distribution that depends only on $\beta$ (but neither on the details of the heat bath nor on the basis $\{b_j\}$). As conjectured in \cite{GLTZ06} and proven using Theorem 1 in \cite{GLMTZ07}, this distribution, the thermal equilibrium distribution of the conditional wave function, is the Gaussian-adjusted-projected (GAP) measure associated with the canonical density matrix of temperature $1/\beta$,
\begin{equation}
\rho_\beta = \frac{1}{Z} e^{-\beta H}\,,\quad Z = \Tr\, e^{-\beta H}\,.
\end{equation}
For any density matrix $\rho$, the measure $GAP(\rho)$ is defined as follows. Let $G(\rho)$ be the Gaussian measure on Hilbert space with covariance $\rho$; multiply $G(\rho)$ by the density function $\|\cdot\|^2$ (adjustment factor) to obtain the measure $GA(\rho)$; project $GA(\rho)$ to the unit sphere in Hilbert space to obtain $GAP(\rho)$.

We now turn to the proof of Theorem~1.

\section{First Part of the Proof: Construction of $U_{ij}$}

We write $M_j$ for the $j$-th column of any $n\times n$ matrix $(M_{ij})$ and
\begin{equation}
  \scp{M_j}{M_\ell} = \sum_{i=1}^n M_{ij}^* M_{i\ell}\,, \quad
  \|M_j\|^2 = \sum_{i=1}^n |M_{ij}|^2\,.
\end{equation}

For $i,j=1,\ldots,n$ let $G_{ij}$ be i.i.d.\ complex Gaussian random
variables with mean 0 and variance 1. To the $n$ columns of the
matrix $(G_{ij})$ apply the Gram--Schmidt orthonormalization
procedure, and call the resulting matrix $(U_{ij})$. That is,
\begin{equation}
  U_{ij} = \frac{G_{ij} - \Delta_{ij}}{\|G_j - \Delta_j\|}
\end{equation}
with
\begin{equation}
  \Delta_{ij} = \sum_{\ell=1}^{j-1} \scp{G_j}{U_\ell} U_{i\ell}\,.
\end{equation}
The procedure
fails if the columns of $(G_{ij})$ are linearly dependent, but
this event has probability 0. Then, as also remarked in \cite{PR04}, $(U_{ij})$ is Haar($U(n)$)
distributed because its first column is uniformly distributed over
the unit sphere in $\CCC^n$, the distribution of the second column
conditional on the first column is uniform over the unit sphere in
the orthogonal complement of the first column, ..., the distribution
of the $j+1$-st column conditional on the first $j$ columns is
uniform over the unit sphere in the orthogonal complement of the
first $j$ columns---and this is exactly the Haar measure.

Our method of proof is to show that $|\sqrt{n}U_{ij}-G_{ij}|$ is in
fact small if $n$ is large. More precisely, we show that for every
$\varepsilon>0$,
\begin{equation}\label{convprob}
  \PPP\left( \sum_{i,j=1}^k |\sqrt{n}U_{ij}-G_{ij}| < \varepsilon \right) \to 1
\end{equation}
as $n\to\infty$. This is called convergence in probability, and to
obtain the claim of the theorem we use the known fact
\cite[Theorem~25.2, p.~284]{Bill} that convergence in probability
implies weak convergence (of the joint distribution of
$\sqrt{n}U_{ij}$ for $i,j=1,\ldots,k$), provided that all random
variables are defined on the same probability space. Here, we can
assume that for all $i,j\in \NNN$, the $G_{ij}$ are defined on the
same probability space.

\section{Second Part of the Proof:\\ Probable Geometry}

The proof of \eqref{convprob} is based on the following
observations:
\begin{itemize}
\item Any two different columns of $(G_{ij})$ tend to be nearly orthogonal.
\item Every column of $(G_{ij})$ tends to have norm close to $\sqrt{n}$.
\item The size of every single entry, $|G_{ij}|$, stays bounded as $n$ grows.
\end{itemize}
These statements are to be understood in the sense that they are fulfilled with high probability for sufficiently large $n$. We now make them precise.

Fix a (small) $\delta>0$. Choose $R>0$ so large that
\begin{equation}\label{defR}
  \PPP \bigl( |G_{ij}|<R \bigr) \geq 1-\delta\,.
\end{equation}
Define the following events corresponding to the three bullets
above:
\begin{equation}
  A_{j\ell}^n := \left\{ \bigl| \scp{G_j}{G_\ell} \bigr| < \sqrt{\frac{n}{\delta}} \right\}
\end{equation}
\begin{equation}
  B_j^n := \left\{  \Bigl| \frac{\|G_j\|^2}{n} -1 \Bigr| < \sqrt{\frac{2}{n\delta}} \right\}
\end{equation}
\begin{equation}
  C_{ij}^n := \left\{ |G_{ij}|<R \right\}
\end{equation}
for $i,j,\ell\leq k$. ($C_{ij}^n$ actually does not depend on $n$, but
never mind.) Each of these events has at least probability
$1-\delta$: $A_{j\ell}^n$ and $B_j^n$ by Chebyshev's inequality and
$C_{ij}^n$ by \eqref{defR}. Thus, the event
\begin{equation}
  D^n:= \bigcap_{\substack{j,\ell=1\\j\neq\ell}}^k A_{i\ell}^n \cap \bigcap_{j=1}^k B_j^n \cap \bigcap_{i,j=1}^k C_{ij}^n
\end{equation}
has at least probability $1-2k^2\delta$, as $2k^2$ is the number of
intersecting sets.

We now show that for sufficiently large $n$, $D^n \subseteq E^n$, where
\begin{equation}
E^n:=\left\{ \sum_{i,j=1}^k |\sqrt{n}U_{ij}-G_{ij}| < \varepsilon \right\}
\end{equation}
is the event in the brackets of \eqref{convprob}. Since
$\delta$ was arbitrary, this fact implies \eqref{convprob}. The
remainder of the proof makes no reference to probabilities, but
concerns only the inclusion $D^n\subseteq E^n$, which can be
regarded as an inclusion between subsets of $\CCC^{n^2}$. Also
$A_{j\ell}^n$, $B_j^n$, and $C_{ij}^n$ will from now on be
regarded as subsets of $\CCC^{n^2}$. (Now the upper index $n$ in the
notation $C_{ij}^n$ becomes useful.) We thus
regard $G_{ij}$ as fixed numbers, and assume that the matrix $(G_{ij})$ lies in the set $D^n$.
When we refer to ``the condition $A_{j\ell}^n$'' we mean the condition that
the $n \times n$ matrix $(G_{im})$ lies in the set $A_{j\ell}^n$.

We proceed to show, by induction over $j\in\{1,\ldots,k\}$, that for
sufficiently large $n$ we have that for all $(G_{ij})\in D^n$, and for all $i=1,\ldots,k$,
\begin{equation}\label{IA2}
  \bigl| \sqrt{n} U_{ij} - G_{ij} \bigr| < \frac{\varepsilon}{k^2}
\end{equation}
and there are constants $C_1, \ldots, C_k>0$ such that for
sufficiently large $n$
\begin{equation}\label{IA1}
  \bigl\| \sqrt{n} U_j - G_j \bigr\| < C_j\,.
\end{equation}
From \eqref{IA2} we see that $(G_{ij})\in E^n$, which is what we need to show.
This induction is the contents of the next, and last, section.

\section{Third Part of the Proof: Estimates}

For $j=1$, note that $U_1 = G_1/\|G_1\|$. By conditions $B_1^n$ and
$C_{i1}^n$,
\begin{equation}
  \bigl| \sqrt{n} U_{i1} -G_{i1} \bigr| =
  \Bigl| \frac{\sqrt{n}}{\|G_1\|} -1 \Bigr| \, |G_{i1}|
  < \frac{2}{\sqrt{\delta n}} R < \frac{\varepsilon}{k^2}
\end{equation}
for sufficiently large $n$. By condition $B_1^n$,
\begin{equation}
  \bigl\| \sqrt{n} U_1 - G_1 \bigr\| =
  \Bigl| \frac{\sqrt{n}}{\|G_1\|} -1 \Bigr| \, \|G_{1}\|
  < \frac{2}{\sqrt{\delta n}} 2 \sqrt{n} = \frac{4}{\sqrt{\delta}} =: C_1\,.
\end{equation}

We now collect four estimates. For $\ell<j$ and sufficiently large
$n$ we find
\begin{equation}
  \bigl| \scp{G_j}{\sqrt{n}U_\ell} \bigr| \leq
  \Bigl| \scp{G_j}{\sqrt{n}U_\ell -G_\ell} \Bigr| +
  \Bigl| \scp{G_j}{G_\ell} \Bigr| \leq
\end{equation}
\begin{equation}\label{est1}
  \leq \|G_j\| \, \|\sqrt{n} U_\ell - G_\ell\| + \sqrt{n/\delta}
  < 2\sqrt{n} C_\ell+ \sqrt{n/\delta} =: C'_\ell \sqrt{n}
\end{equation}
where we have used the Cauchy--Schwarz inequality, $A_{j\ell}^n$,
$B_j^n$, and the induction hypothesis \eqref{IA1}. As the next
estimate, for $i\leq k$,
\begin{equation}
  |\Delta_{ij}| \leq \sum_{\ell=1}^{j-1} \frac{\bigl| \scp{G_j}{\sqrt{n}U_\ell} \bigr|}
  {\sqrt{n}} \frac{|\sqrt{n} U_{i\ell}|}{\sqrt{n}} \leq
\end{equation}
\begin{equation}
  \leq \sum_{\ell=1}^{j-1} C'_\ell \frac{1}{\sqrt{n}}
  \Bigl( |\sqrt{n} U_{i\ell} - G_{i\ell}| + |G_{i\ell}| \Bigr) <
\end{equation}
\begin{equation}\label{est2}
  < \Bigl(\sum_{\ell=1}^{j-1} C'_\ell\Bigr)
  \Bigl( \frac{\varepsilon}{k^2} + R \Bigr) \frac{1}{\sqrt{n}}
  =: \frac{C''_j}{\sqrt{n}}
\end{equation}
using \eqref{est1}, the induction hypothesis \eqref{IA2}, and
$C_{i\ell}^n$. As the third estimate, for $j \leq k$
\begin{equation}
  \|\Delta_j\|^2 = \sum_{i=1}^n |\Delta_{ij}|^2 \leq \frac{1}{n}
  \sum_{i=1}^n \sum_{\ell=1}^{j-1} \bigl| \scp{G_j}{\sqrt{n}U_\ell} \bigr|^2 \,
  |U_{i\ell}|^2 <
\end{equation}
\begin{equation}\label{est3}
  < \frac{1}{n}\sum_{\ell=1}^{j-1} (C'_\ell)^2 n \|U_\ell\|^2 =
  \sum_{\ell=1}^{j-1} (C'_\ell)^2 =: C_j'''
\end{equation}
using \eqref{est1} and the fact that $U_\ell$ is a unit vector. As
the last estimate,
\begin{equation}\label{est4}
 \left| \frac{\sqrt{n}}{\|G_j-\Delta_j\|} -1 \right| <
 \left( 2\sqrt{\frac{2}{\delta}} + 2 \sqrt{C_j'''} \right) \frac{1}{\sqrt{n}}
 =: \frac{\tilde{C}_j}{\sqrt{n}}\,,
\end{equation}
which is easily obtained from $\|G_j\| - \|\Delta_j\| \leq \|G_j -
\Delta_j\| \leq \|G_j\| + \|\Delta_j\|$ and \eqref{est3} in the
following way:
\begin{equation}
  \frac{\sqrt{n}}{\|G_j-\Delta_j\|} - 1 \leq \frac{\sqrt{n}}{\|G_j\| - \|\Delta_j\|} - 1 <
\end{equation}
\begin{equation}
  <\frac{\sqrt{n}}{\|G_j\| - \sqrt{C_j'''}} - 1= \frac{1}{\|G_j\|/\sqrt{n} - \sqrt{C_j'''/n}} -1\,.
\end{equation}
Since, using $B_j^n$, $\|G_j\|/\sqrt{n} > \sqrt{1-\sqrt{2/n\delta}} >
1- \sqrt{2/n\delta}$, and since
\begin{equation}
  \frac{1}{1-x} < 1+2x
\end{equation}
for sufficiently small $x>0$, we obtain that, for sufficiently large
$n$,
\begin{equation}
  \frac{\sqrt{n}}{\|G_j-\Delta_j\|} - 1< \frac{1}{1- \sqrt{2/n\delta} - \sqrt{C_j'''/n}} -1
  < 2\sqrt{\frac{2}{n\delta}} + 2 \sqrt{\frac{C_j'''}{n}}\,.
\end{equation}
Together with an (even narrower) lower bound obtained by similar
arguments, this yields \eqref{est4}.

From these four estimates, the first induction claim \eqref{IA2}
follows for $j\leq k$ because, for sufficiently large $n$,
\begin{equation}
  \bigl| \sqrt{n} U_{ij} - G_{ij} \bigr| =
  \Bigl| \frac{\sqrt{n}}{\|G_j -\Delta_j\|} (G_{ij} - \Delta_{ij}) - G_{ij} \Bigr| \leq
\end{equation}
\begin{equation}
  \leq \Bigl| \frac{\sqrt{n}}{\|G_j - \Delta_j\|} -1 \Bigr| \, |G_{ij}| +
  \frac{\sqrt{n}}{\|G_j - \Delta_j\|} |\Delta_{ij}| <
\end{equation}
\begin{equation}
  < \frac{\tilde{C}_j}{\sqrt{n}} \, R + 2 \frac{C''_j}{\sqrt{n}} < \frac{\varepsilon}{k^2}\,,
\end{equation}
where we have used \eqref{est4}, $C_{ij}^n$, \eqref{est4} with
$\tilde{C}_j/\sqrt{n} <1$, and \eqref{est2}. The second induction
claim \eqref{IA1} follows from
\begin{equation}
  \bigl\|\sqrt{n}U_{j} - G_{j}\bigr\| \leq \Bigl| \frac{\sqrt{n}}{\|G_j -\Delta_j\|} -1 \Bigr|
  \, \|G_j\| +   \frac{\sqrt{n}}{\|G_j -\Delta_j\|} \|\Delta_j\| <
\end{equation}
\begin{equation}
  < \frac{\tilde{C}_j}{\sqrt{n}} 2\sqrt{n} + 2 \sqrt{C_j'''} = 2\tilde{C}_j + 2 \sqrt{C_j'''}
  =: C_j\,.
\end{equation}
This completes the proof.


\bigskip

We close with a remark on the parenthesis in Theorem 1: ``the upper left (or, in fact, any) $k\times k$ submatrix.'' We elucidate the meaning of ``any.'' To select a $k\times k$ submatrix means to select $k$ rows and $k$ columns. This selection must be deterministic (\y{i.e., non-random,} or at least independent of the $U_{ij}$) but may depend on $n$. \y{Indeed,} if the selection depended on the $U_{ij}$, one could, for example, select those rows and columns where $U_{ij}$ happens to be exceptionally close to zero, which would lead to a different asymptotic distribution. On the other hand, for a selection depending on $n$, Theorem 1 remains true: \y{to see this, recall that for a compact group such as $U(n)$, the Haar measure is both left-invariant and right-invariant; as a consequence, Haar$(U(n))$ is invariant under any (non-random) permutation of either the rows or the columns, and thus} all $k\times k$ submatrices have the same distribution.

\end{document}